\documentclass[1p]{elsarticle}

\usepackage{lineno}
\modulolinenumbers[5]
\pdfoutput=1

\makeatletter
\def\ps@pprintTitle{%
   \let\@oddhead\@empty
   \let\@evenhead\@empty
   \def\@oddfoot{\reset@font\hfil\thepage\hfil}
   \let\@evenfoot\@oddfoot
}
\makeatother
\journal{European Journal of Operations Research}

\biboptions{authoryear}

\usepackage{amssymb}
\usepackage{amsmath}
\usepackage{algorithm,algorithmic}
\DeclareMathOperator*{\argmax}{argmax}
\usepackage{multirow}

\usepackage{subfigure}

\usepackage{url}

\usepackage{array,tabularx}

\begin{document}

\begin{frontmatter}

\title{E-commerce warehousing: learning a storage policy}

\author[label1,label2]{Adrien Rim\'el\'e\corref{cor1}}
\cortext[cor1]{Corresponding author}
\ead{adrien.rimele@polymtl.ca}

\author[label5]{Philippe Grangier}
\author[label1,label3]{Michel Gamache}
\author[label1,label2]{Michel Gendreau}
\author[label1,label2]{Louis-Martin Rousseau}

\address[label1]{Department of Mathematics and Industrial Engineering, Polytechnique Montreal}
\address[label2]{CIRRELT, Research Centre on Enterprise Networks, Logistics and Transportation}
\address[label3]{GERAD, Group for Research in Decision Analysis}
\address[label5]{IVADO Labs}

\begin{abstract}
E-commerce with major online retailers is changing the way people consume. The goal of increasing delivery speed while remaining cost-effective poses significant new challenges for supply chains as they race to satisfy the growing and fast-changing demand. In this paper, we consider a warehouse with a Robotic Mobile Fulfillment System (RMFS), in which a fleet of robots stores and retrieves shelves of items and brings them to human pickers. To adapt to changing demand, uncertainty, and differentiated service (e.g., prime vs. regular), one can dynamically modify the storage allocation of a shelf. The objective is to define a dynamic storage policy to minimise the average cycle time used by the robots to fulfil requests. We propose formulating this system as a Partially Observable Markov Decision Process, and using a Deep Q-learning agent from Reinforcement Learning, to learn an efficient real-time storage policy that leverages repeated experiences and insightful forecasts using simulations. Additionally, we develop a rollout strategy to enhance our method by leveraging more information available at a given time step. Using simulations to compare our method to traditional storage rules used in the industry showed preliminary results up to 14\% better in terms of travelling times.
\end{abstract}


\end{frontmatter}

%\linenumbers

\section{Introduction}
\label{introduction}

\subsection{Warehousing in e-commerce}
\label{warehousing}

Warehousing occupies a central role in a supply chain. According to \cite{Gu2007}, the primary purposes of a warehouse are to act as a buffer to adapt to the variability of the production flow; to consolidate products; and to add marginal value such as pricing, labelling or customisation. 

The recent and massive development of e-commerce introduces great challenges in the \emph{Business to Customer} segment. E-commerce involves enormous volumes of orders, and online sales keep growing in number: in 2016, e-commerce sales grew by 23.7\%, accounting for 8.7\% of the total market \citep{Boysen2019ManualProblem}. These orders consist, on average, of very few items: the same authors mention that the average number of items per order at Amazon is only 1.6. E-commerce also faces great demand uncertainty, with fast-changing demand trends, and the need to satisfy orders quickly \citep{Yaman2012ReleaseDelivery, Boysen2019ManualProblem, Boysen2019WarehousingSurvey}. Amazon, for instance, offers customers \emph{Prime}, a differentiated service that used to offer delivery within two days and is now, under certain conditions, promising same-day delivery. Moreover, with easily accessible alternatives available to customers, online retailers face high competitive pressure, and the link between logistics performance and customer loyalty is much stronger compared to other industries \citep{Ramanathan2010TheE-commerce}.

Parts-to-picker Automated Storage and Retrieval Systems (AS/RS), with the typical single, aisle-captive crane retrieving bins from a static rack, have been used for several decades and have proven their superior operational efficiency compared to manual systems, at the cost of an essential initial investment \citep{Roodbergen2009}. However, some limitations of AS/RS become apparent when facing new challenges posed by the rise of e-commerce: cranes are sensitive points of failure; expandability is complicated and expensive; batching and zoning require consolidation and a delay-sensitive dependency between operators; and a single crane per aisle can only perform so many cycles. According to \cite{Davarzani2015}, ``adaptation is urgent'' to keep up with the growing demand.

\subsection{Robotic Mobile Fulfillment System}
\label{RMFS}

Kiva Systems created the first RMFS in 2006, this company was later bought by Amazon and re-branded Amazon Robotics in 2012 \citep{Wurman2008CoordinatingWarehouses, Azadeh2017RobotizedOpportunities,Boysen2019WarehousingSurvey}. However, RMFSs, also referred to as rack-moving mobile robot based warehouses, or Kiva warehouses, are not limited to Amazon warehouses. Other companies have developed their own systems: for instance, Alibaba's Zhu Que robots, the Quicktron robots (Huawei), the Open Shuttles by Knapp, CarryPick by Swisslog, Butler by GreyOrange, the Locus Robotics system and others \citep{Kirks2012CellularSystems, Boysen2019WarehousingSurvey}.

This new type of automated warehouse involves a large fleet of small robots that move freely around the storage area to retrieve shelves full of items and bring them to operators for picking (possibly after waiting in a queue) before returning to a storage location. Interestingly, when a robot is loaded, it must travel along the aisles, but when unloaded, it can pass under stored shelves to take shortcuts. This type of operations is equivalent to a mini-load with dual-command cycles system with a fleet of non-aisle-captive robots \citep{Roodbergen2009}. For a more detailed description of an RMFS and its advantages compared to traditional systems, the interested reader is referred to the following references: \cite{Wurman2008CoordinatingWarehouses}; \cite{DAndrea2008FutureFacilities}; \cite{Enright2011OptimizationSystems}; and \cite{Azadeh2017RobotizedOpportunities}.

\cite{Rimele2020RoboticApplications} present a mathematical modelling framework that formalises optimisation opportunities in an RMFS. Because e-commerce promises fast delivery of small orders, the fulfilment of requests follows what can be called an \emph{order streaming} logic \citep{Banker2018NewAutomation}, which ``drops orders [...] to the floor as soon as they are received''. Because of this specificity, the option of batching requests is minimal, and the system needs to account for uncertainty and adapt online to newly revealed orders. Among the operational decisions, the real-time storage allocation problem can leverage the great storage flexibility of the system. When a robot leaves a picking station, it does not need to store the shelf in its initial location but can choose \emph{any} free location. This allows for the real-time adaptation of the storage layout in response to the demand, thereby minimising the travelling time of the robots, which translates into increasing an upper bound of the throughput.

\subsection{Contributions and paper structure}
\label{contribution}

This work focuses on a specific operational control aspect of a Robotic Mobile Fulfillment System for e-commerce, namely the problem of learning an effective dynamic storage policy. When a robot is available at an operator's station, a decision must be made in real-time about which storage location to use to minimise the average cycle time of the robots. To make such a decision, one can use information about the next request to process, statistical insights about the demand, and other input regarding the current state of the warehouse's layout. This information at a given time step defines a state's representation in an environment that we approximate as a partially observable Markov decision process (POMDP). We propose the application of a variant of a Q-learning agent \citep{Watkins1989LearningRewards,Hessel2018Rainbow:Learning, Mnih2015}, a Reinforcement Learning (RL) algorithm, to learn an efficient dynamic storage policy that aims at minimising the average cycle time. More precisely, we implement a Differential Double Q-learning agent with experience relays to minimise the average travelling time of the robots. Inspired by \emph{Class-based} storage policies, the decision on where to return a shelf is made based on a discretisation of the storage area into zones. A discrete event simulator is used to evaluate the performance of our method compared to baselines on different levels of skewed demand distributions. To further improve the performance of our method, we propose an extension using a rollout strategy over a fixed horizon of already-revealed requests. This strategy is also applicable to the baseline methods from the literature.

The presentation of our proposed approach is organised as follows. Section \ref{problem_definition} first formalises the optimisation problem of dynamically allocating storage locations to shelves, then presents a literature review related to storage policies in automated warehouses and operational control applications specific to RMFS. Section \ref{methodology} presents our methodology, with a discretisation of the storage area into zones, the system representation as a POMDP, and the reinforcement learning variant of a Q-learning agent used to learn an efficient storage policy, as well as the rollout strategy to further improve performance. Section \ref{sec:simulation_study} explains how we generated data in our simulations and shows results obtained by our method compared with baselines from the literature. Finally, Section \ref{conclusions} summarises the obtained results and discusses future research avenues.

\section{Problem definition and literature}
\label{problem_definition}

\subsection{Decision-making framework}
\label{decision_making_framework}

In this section, we present the storage allocation decision-making framework, as well as the relevant literature on the topic.

\begin{figure}
\centering
\includegraphics[width=110mm]{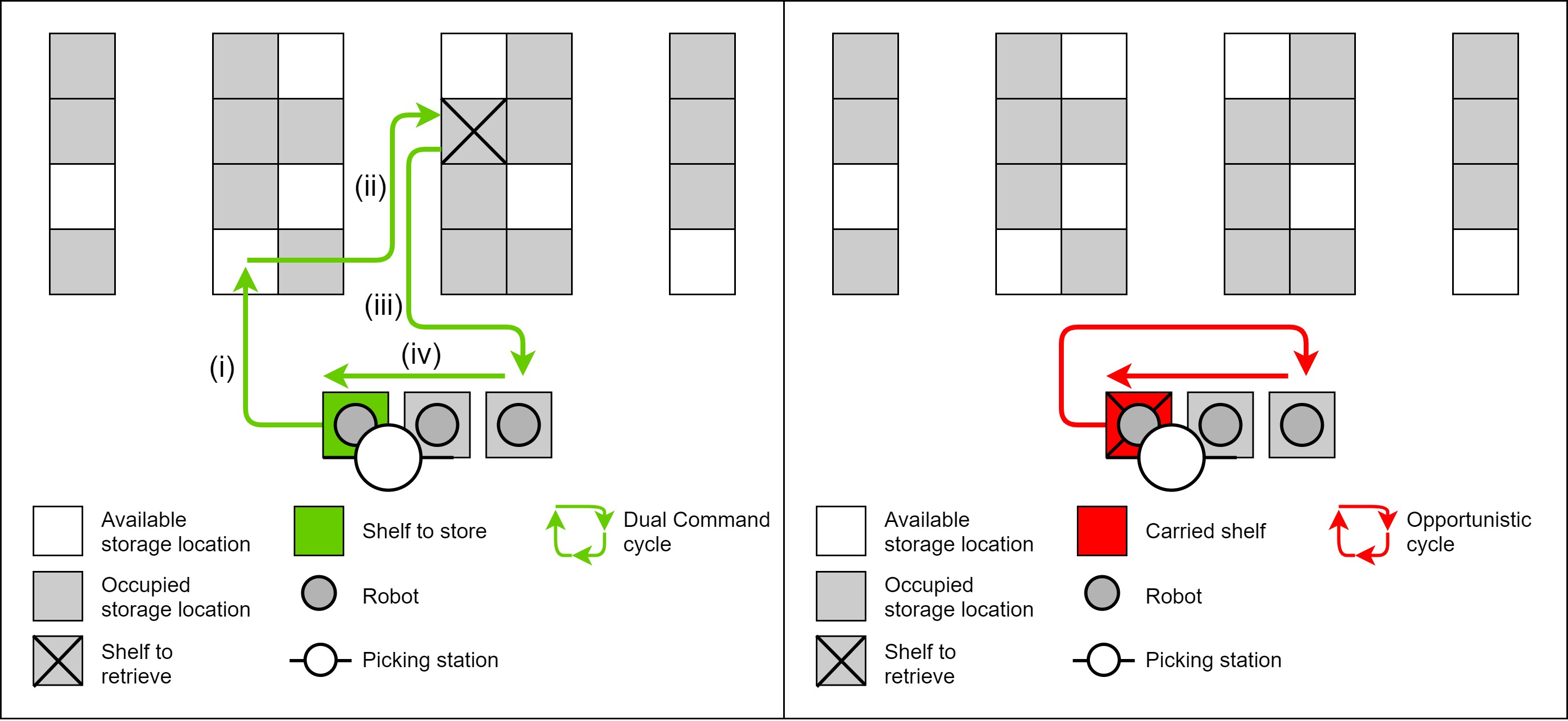}
\caption{Dual-command cycle (left) and opportunistic task (right)}
\label{fig:cycles}
\end{figure}

For the case of an e-commerce Amazon-type warehouse, a modelling of operational decisions is proposed by \cite{Rimele2020RoboticApplications}, with a stochastic dynamic optimisation formulation. In this formulation, an incoming request corresponds to a single item type and consolidation of orders, if any, is assumed to take place in the downstream sections of the warehouse. Orders are revealed online and are considered available for fulfilment as soon as they arrive, with no batching. An \emph{event} corresponds to a robot being available for a task. The tasks of a robot correspond to a sequence of \emph{dual-command (DC) cycles} for which the total time can be decomposed into: (i) storage access time; (ii) interleaving (or transition) time to retrieval location; (iii) retrieval time; and (iv) waiting time at the picking station's queue (see arcs on the left side of Figure \ref{fig:cycles}). To create those cycles, storage and retrieval decisions must be made. After the picking operation, the robot is available at the picking station, and it has to perform a \emph{storage task}, which consists in selecting an available location for storage. When the shelf has been stored, the robot needs to perform a \emph{retrieval task}, which consists in choosing which order to fulfil next, from which shelf, and which picking station. Besides those two typical tasks, we also consider an \emph{opportunistic task}, a task that can only occur in a particular situation where a requested item is stored on a shelf currently being carried by the robot. In this case, the robot can go directly to the back of the waiting line at the operator, without passing by the storage area (right side of Figure \ref{fig:cycles}).

In this work, we focus on optimising storage decisions only. It is assumed that an external system provides the sequence of orders and the decisions on which shelf and picking station to use. To be more precise, this system simply considers orders already revealed and yet to be fulfilled and dynamically sequences them by increasing deadlines, as an emergency rule. If an order requires a shelf carried by another robot, this order is temporarily skipped to be processed by that other robot as an opportunistic task. Also, opportunistic tasks are performed as often as possible. The objective is to minimise the average travelling time, which is directly related to the maximal throughput capacity of the warehouse \citep{Park2012OrderModels}. We make other simplification assumptions similar to the ones enumerated in the simulation study in \cite{Rimele2020RoboticApplications}. We consider one picking station only, an item is associated with exactly one shelf in the considered storage area, replenishment is explicitly ignored even if some orders could correspond to replenishment tasks, travelling and processing times are deterministic and, finally, congestion of robots in the aisles is ignored.

\subsection{Storage policies}
\label{storage_policies}

Commonly used and studied storage policies for AS/RS are \emph{Random} storage policy or policies based on the container's turnover rate. The former gives a container an equal probability to occupy any available location. While this method is space-efficient and easy to implement in practice, it can obviously result in poor accessibility and low productivity during operations. The latter method is typically expressed as a full-turnover-based storage assignment or a class-based storage assignment. A full-turnover-based assignment locates the most-requested containers closest to the Input/Output (I/O) point; however, a practical concern appears when new containers enter and leave the system, or when the turnover rates change. In \emph{Class-based} storage, storage locations are clustered into a given number of classes based on their proximity to the I/O point (see left side of Figure \ref{fig:class_zone_layouts}). Then, a container is assigned to a zone based on its turnover rate; its storage within this zone is random. While this method benefits from both \emph{Random} storage and full-turnover-based assignments, it requires the definition of the number of classes and their dimensions. Numerous analytical and simulation studies such as \cite{Hausman1976OptimalSystems}, \cite{Graves1977Storage-RetrievalSystems}, \cite{Schwarz1978SchedulingResults}, \cite{Bozer1984Travel-TimeSystems}, and \cite{Gagliardi2012OnSystems} show that \emph{Class-based} storage assignment demonstrates the best performance and practicality. However, while those applications generally apply to both single and dual-command cycles system, it seems they have been more driven by the former. In another simulation study, \cite{VandenBerg2000SimulationSystem} find out that for dual-command cycles, selecting the storage location that minimises the cumulative distance from the I/O point to the storage location, plus the distance from the storage location to the following retrieval location, gives the best average completion time. This decision rule is often referred to as the \emph{Shortest Leg} storage policy.

\subsection{RMFS-related work}
\label{rmfs_related_work}
Some recent work focuses more specifically on RMFS, either on storage decisions or closely related topics like performance measures and others. 
\cite{Lamballais2017EstimatingSystem} use queueing network models to analytically estimate maximum order throughput, average cycle times, and robot utilisation. Using these models allows one to quickly evaluate different layouts or robot zoning strategies. Among their results, they show the good accuracy of their models compared to simulations. They also show that the throughput is quite insensitive to the length-to-width ratio of the storage area, but is more strongly affected by the location of the picking stations. Also using queueing network models, \cite{Roy2019} study the effect of dedicating robots exclusively to either retrieval tasks or replenishment tasks, or a combination of both. They find that the latter reduces the picking time by up to one third, but increases the replenishment time up to three times. They also study the allocation to multiple classes and conclude that allocating robots to least congested classes results in a similar performance as a dedicated zone policy.

\cite{Boysen2017Parts-to-pickerEnvironment} study the problem of synchronising the processing of orders at a picking station with the visits of the shelves carried by robots. By considering a given set of orders, and a capacity to process orders simultaneously, their objective is to minimise the number of visits of shelves to the picking station. They propose a MIP model, a decomposition approach, and different heuristics and show that they can halve the robot fleet size.

Regarding storage allocation, \cite{Weidinger2018StorageWarehouses} define a \emph{Rack Assignment Problem}, which assigns each stopover (return of a rack) of a given set (batching) to an open storage location. They consider the time at which each rack visits the picking station to be known. They propose a MIP model as a special case of an interval scheduling problem, with the surrogate objective of minimising the total loaded distance. The model formalises the fact that two stopovers can occupy the same storage location only if their storage intervals cannot overlap. They do not consider robots individually, but instead they expect that a robot will always be available for the task. They propose solving their model using an Adaptive Large Matheuristic Search (ALMS), and compare their results with other storage policies, such as those presented in Section \ref{storage_policies}. On their surrogate objective, ALMS demonstrates good performance, outperforming the other policies on the total travel time (unloaded robots are considered to travel twice as fast as loaded ones). Similar to the finding of \cite{VandenBerg2000SimulationSystem}, the authors note the very good performance of the \emph{Shortest Leg} policy; compared to their method, it only increases the travel distance by 3.49\% and the size of the robot fleet by 2.17\%, without batching.

\cite{Yuan2018} propose a fluid model to assess the performance of velocity-based storage policies in terms of travelling time to store and retrieve shelves. They find that when the items are stowed randomly within a shelf, a 2-class or 3-class storage policy reduces the travelling time between 6 to 12\%. This theoretical result is validated by simulation. Additionally, they find that if items are stowed together based on their velocity, the travelling distances considered can be reduced by up to 40\%.

In situations when batching or zoning requires consolidation, different options can be considered. \cite{Boysen2019ManualProblem} consider a manual consolidation with put walls (sorting shelves). They propose a MIP to optimise the release sequence of bins from intermediate storage so that orders are quickly sorted on the put wall and idle times of operators are minimised. In a case where the intermediate storage uses an AS/RS that serves the consolidation area with a conveyor belt, \cite{Boysen2018OptimizingProblem} propose a MIP model for the minimum order spread sequencing problem, which consists in minimising the number of conveyor segments between the first and last occurrence of an order.

\cite{Merschformann2019DecisionSystems} use a discrete event simulator presented in \cite{Merschformann2018RAWSim-O:Systems}, which considers most, if not all, of the operational decisions occurring in an RMFS. They test a multitude of combinations of decision rules concerning, for instance, the assignment to a picking or replenishment station, the robot selection, storage allocation, etc. They observe significant performance differences between distinct combinations, as well as strong cross-dependencies between some decision rules, suggesting there may be some benefit to exploring integrated approaches.

\section{Methodology}
\label{methodology}

This section presents our method for solving the dynamic storage allocation problem in an RMFS. We first define the notion of zones that we use in our storage policy. Then, we explain the representation of the system as a POMDP and the features extraction. In Section \ref{sec:rl_method}, we describe the variant of Q-learning, a Reinforcement Learning agent, that we use. Section \ref{sec:imple_details} presents implementation details essential to the success of the method. Finally, Section \ref{sec:rollout} describes a look-ahead rollout strategy to enhance our storage policy, by making better-informed decisions using already-revealed orders.

\subsection{Storage zones definition}
\label{zones_definition}

\begin{figure}
\centering
\subfigure{\resizebox*{5.5cm}{!}{\includegraphics{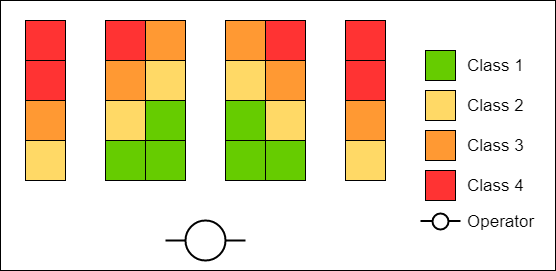}}}\hspace{5pt}
\subfigure{\resizebox*{5.5cm}{!}{\includegraphics{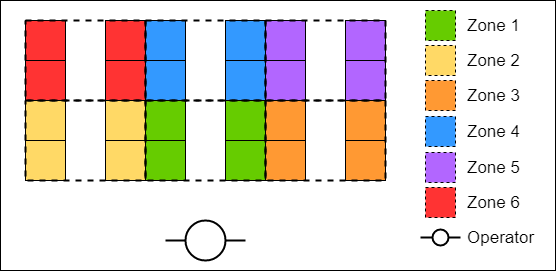}}}
\caption{Concentric class-based storage layout (left) vs. proposed zones layout (right)} \label{fig:class_zone_layouts}
\end{figure}

In this work, inspired by \emph{Class-based} storage policies from the literature, we propose optimising a policy that will allocate shelves to classes instead of exact locations. Contrary to the standard \emph{Class-based} storage assignment presented in Section \ref{storage_policies}, the online assignment to a class is not only made based on the turnover rate of the shelf to store, but also on a set of diverse features. Since we are interested in minimising dual-command cycle completion times, as opposed to the access to storage or retrieval time alone (like in single command cycles), the interleaving time between a storage location and the next retrieval location is essential. With the usual definition of classes distributed with an increasing distance from the picking station, such as in Figure \ref{fig:class_zone_layouts} (left), the access and retrieval times depend directly on the assigned class, but the interleaving time is quite random. For instance, in the ideal case where we store and retrieve shelves from the same storage class 2, we see that the cycle can either be very good (if the two locations are adjacent), or very bad if the locations are at the two extremities. To better account for the importance of interleaving time, we propose another definition of classes, which we now call \emph{zones} as represented in Figure \ref{fig:class_zone_layouts} (right). These zones are now clustered into regular shapes such that locations within the same zone are close to each other and such that the zones are distributed in an increasing distance from the picking station. With this new zone definition, a storage and a retrieval from the same zone would entail a short interleaving time and, thus, a good cycle. Note that the selection of the exact storage location within a zone is arbitrarily made by choosing the open location closest to the picking station.

\subsection{System representation}
\label{system_representation}

\cite{Rimele2020RoboticApplications} model the storage and retrieval operations in an AS/RS as a Markov Decision Process (MDP), which gives a full representation of the warehouse, including all of the revealed orders, as a state of the system when an event occurs (corresponding to the need for a storage decision here). Instead of giving such a complete representation of the current state of the warehouse, we propose a \emph{Partially Observable} MDP (POMDP), where only some characteristics of the real state are used as state representation. In the traditional \emph{Class-based} storage policy, for instance, the representation of the state is simply the turnover rate of the shelf to store, and based on this information the robot knows which class to assign it to.

We decided to represent a state with some more characteristics (features), given below:
\begin{itemize}
\setlength\itemsep{0.2em}
\item Average turnover rate of the shelf to store
\item Relative rank of the shelf's turnover rate
\item Zone of the next retrieval task (one-hot encoding, one feature per zone)
\item Occupation levels of each zone (one feature per zone)
\item For each zone, the number of robots that are currently moving toward the zone to retrieve a shelf
\end{itemize}

An extra feature is used to encode whether the next task is an opportunistic one (the shelf is not stored but is already being carried by the robot). In total, if $n$ denotes the number of zones, there are $3n+3$ scalars in the features vector.

The intuitive idea is that there must be a balance between greedily minimising the immediate cycle time, similar to the \emph{Shortest Leg} storage policy, and minimising the future access time, similar to a turnover-based storage policy as presented in Section \ref{storage_policies}. The objective is to learn a policy that, based on the above set of features, will aim at minimising the average cycle time.

\subsection{Reinforcement learning method}
\label{sec:rl_method}

We propose using a Reinforcement Learning (RL) method to learn an efficient storage policy. The field of RL offers a diversity of methods to solve decision problems whose objectives are defined recursively by a Bellman equation, introduced in dynamic programming (DP). Compared to traditional DP methods that require complete information about the model and probability distributions, RL methods can learn from repeated experiences, even in large state and action space, when using function approximations. As presented by \cite{Sutton2018ReinforcementIntroduction}, in RL the decision-maker, or agent, interacts with its environment. At every time step $t$, and based on the current representation of the state $S_t$, the agent selects an action $A_t$, which impacts the environment and receives from it a reward $R_{t+1}$ and a new state representation $S_{t+1}$. Based on repeated experiences, the goal of the agent is to learn a policy that maximises the onward cumulative sum of rewards. Note that the agent is not necessarily given the complete state characteristics; it can only be given a partial representation through a POMDP. The state representation defined in Section \ref{system_representation} is the partially observable state given to the agent in our approach. In value-based RL, given a policy $\pi$, the agent learns the State-Action value function $Q_{\pi}(s,a)=\mathbb{E}\left [ R_1 + R_2 + ... \:|\: S_0=s, A_0=a, \pi \right ]$ that represents the expected (possibly discounted) cumulative reward at state $s$ if action $a$ is taken, if policy $\pi$ is followed afterwards. The optimal value function $Q^*$ is such that $Q^*(s,a) = \textup{max}_{\pi}Q_{\pi}(s,a)$ for all state and action pairs, and obtaining a decision policy from a value function simply consists in acting greedily with regard to action values : $A_{\pi}(s)=\textup{argmax}_{a}Q_{\pi}(s,a)$. 

In Q-learning \citep{Watkins1989LearningRewards}, the value function is iteratively updated by sampling based on the Bellman equation: $Q(S_t,A_t) \leftarrow Q(S_t,A_t) + \alpha [ R + \gamma \: \textup{max}_{a}Q(S_{t+1}, a) - Q(S_t, A_t)]$, so the update uses the received reward and the current estimate of the next state value (bootstrap) as its target. Importantly, for convergence to optimality, the Q-learning agent needs to correctly learn the values of its actions (exploitation) while trying new actions, possibly more promising (exploration). This trade-off is central to any RL method, and here we will consider an $\epsilon$-greedy selection rule for this purpose ($\epsilon$ \% of the time, the greedy action is selected, otherwise a random action is taken).

Even if our state space has a small number of dimensions and the number of actions is also small, calculating the value function individually for every state-action pair with a tabular method would be out of reach considering the exponential number of combinations. Instead, Deep Q-learning uses a neural network (NN) as a function approximator that takes the state features as inputs and outputs the estimated values for each possible action \citep{Mnih2015}. Using function approximators allows to escape the curse of dimensionality by leveraging information about similar states and by generalising to unseen state-action pairs using only a limited number of parameters. Let us denote by $\mathbf{w_t}$ the current parameters, weights, of the neural network function approximator and by $Q(s, a, \mathbf{w_t})$ the value of state $s$ - action $a$ estimated by this neural network. Because a neural network is differentiable, we can minimise the square loss error between the current estimate value and the target by taking a step in the opposite direction of the value function's gradient: $\mathbf{w_{t+1}} \leftarrow \mathbf{w_t} + \alpha [R_{t+1} + \gamma \: \textup{max}_{a} Q(S_{t+1}, a, \mathbf{w_t}) - Q(S_t, A_t,\mathbf{w_t})]\nabla Q(S_t, A_t, \mathbf{w_t})$.

In the case of an infinite time horizon system, such as ours where retrieval tasks never end, rewards cannot simply be accumulated because the state values would never converge but rather would diverge to infinity. The most common solution is to apply a discount factor $\gamma$ in the Bellman equations, such that for a given policy $\pi$ we have: $Q_{\pi}(s)=\mathbb{E}_{\pi}[R_{t+1} + \gamma \: Q_{\pi}(S_{t+1}) \:|\: S_t=s]$, with $0 < \gamma < 1$. The justification for using a discount factor is first for the mathematical convenience of convergence of the cumulative finite rewards. Also, in some systems, a discount factor can represent an economical discount rate, or at least the idea that rewards received now are worth more than they would be later. Yet in our case, no incentive should be given for early performance since operations throughout the day have equal importance and no termination occurs. An alternative to discounted rewards is what is called \emph{differential returns} or \emph{average rewards} \citep{Sutton2018ReinforcementIntroduction}. Introduced in Dynamic Programming \citep{BERTSEKAS1976DynamicControl} and later in Reinforcement Learning (R-Learning in \cite{Watkins1989LearningRewards}), the idea is to optimise the average sum of reward $\lim_{n \to \infty}\frac{\sum_{t=1}^{n} R_t}{n}$ by considering return values equal to $R_t -\bar{R}$, with $\bar{R}$ the average reward of policy $\pi$. This setting has the interesting property of having the sum of returns converging (to 0) while maximising the average reward in the process. When action $A_t$ of policy $\pi$ is taken, the state-action value update rule becomes $Q(S_t,A_t) \leftarrow Q(S_t,A_t) + \alpha [ R - \bar{R} + \textup{max}_{a}Q(S_{t+1}, a) - Q(S_t, A_t)]$.

Q-learning is known to suffer from an overestimation of actions' values. Two main reasons explain this phenomenon. First, in its update mechanism, the maximum value of the next state-action pair is used as an approximation for the maximum expected value. Second, both the action selection and value estimation come from the same maximum operator that suffers from estimation errors, and thus favours overestimated values. To answer this problem, \cite{VanHasselt2010DoubleQ-learning} proposes a Double Q-learning method in the tabular case, which was later extended to Deep Q-learning in \cite{vanHasselt2016DeepQ-Learning}. For Deep Q-learning, the general idea is to separate the action selection from the value estimation using two neural networks, the online NN and the target NN. The authors demonstrate that this approach eliminates the positive bias on the action values and results, in general, in better performance. The online NN is updated as before, but the target NN is not. Instead, it receives a copy of the weights of the online NN every $K$ iterations.

Algorithm \ref{alg:algo1} combines the different elements presented above into one agent, denoted as \emph{Double Differential Deep Q-Learning}. Step 1 initialises a starting state, the average reward (arbitrarily set, for instance, to 0) and the weights of both the online and target neural networks, equal to each other. For every iteration, Step 3 selects an action, greedily or not, takes it, and observes the reward and the next state. Step 4 calculates the target value for state $S_t$ - action $A_t$ pair. Step 5 does the gradient descent to update the online network. If the action selected $A_t$ was greedy with respect to the online network, Steps 6 to 8 update the average reward $\bar{R}$. Finally, every $K$ iterations, the target network receives a copy of the online network weights (Steps 9 to 11).

\begin{algorithm}
\begin{algorithmic}[1]
\small
\STATE initialise state $S_{0}, \bar{R}, \mathbf{w_t} = \mathbf{w^{-}_{t}}$
\FOR{$t = 0$ to $\infty$}
    \STATE{take action $A_{t}$ ($\epsilon\text{-greedy}$ wrt $Q(S_{t}, ., \mathbf{w_{t}})$), observe $R_{t+1}$, $S_{t+1}$} 
    \STATE{$Y_{t} \leftarrow R_{t+1} - \bar{R} + Q(S_{t+1},\argmax_a Q(S_{t+1}, a, \mathbf{w_{t}}), \mathbf{w^{-}_{t}})$}
    \STATE{$\mathbf{w_{t+1}} \leftarrow \mathbf{w_t} + \alpha [Y_{t}-Q(S_t, A_t,\mathbf{w_t})]\nabla Q(S_t, A_t, \mathbf{w_t})$}
    \IF{$A_t==\argmax _a Q(S_{t}, a, \mathbf{w_{t}})$} 
        \STATE {$\bar{R}\leftarrow \bar{R} + \beta [Y_{t}-Q(S_t, A_t,\mathbf{w_t})]$}
    \ENDIF
    \IF{$t \: \text{mod} \: K == 0$} 
        \STATE {$\mathbf{w^{-}_{t}} \leftarrow \mathbf{w_t}$}
    \ENDIF
\ENDFOR

\end{algorithmic}
\caption{Double Differential Deep Q-Learning}
\label{alg:algo1}
\end{algorithm}

\subsection{Implementation details}
\label{sec:imple_details}

In our experiments, we use a fully connected feed-forward neural network as a function approximator. While the neural network has one output per storage zone, which represents the expected cumulative cycle times onwards (relative to the average cycle time), if the storage zone is selected, some zones may not be selected if they are currently full. For this reason, when we select an action based on its value out of the neural network, we first need to filter out full zones. 

Another point worth mentioning deals with backpropagation when an opportunistic task is performed. Such a task is enforced when possible, which occurs when the next retrieval zone is not a zone but the picking station (encoded in the features extraction). Even though the agent does not require any action to be taken, it is essential to backpropagate the observed reward in the neural network for all of the actions. The reason for that is bootstrapping. If a non-opportunistic task transitions to a state that will enforce an opportunistic task, the state-action value update in Q-learning uses the observed cycle time as well as the estimated value of the future state, which, in this case, requires an opportunistic task. The neural network must then be accurate for such a state and learn that the values of all the storage zones are to be equal to the value of an opportunistic task.

As in \cite{Mnih2015}, we use experience replay. It consists in storing simulation experiences (vector $S_t, A_t, R_{t+1}, S_{t+1}$) into a fixed size replay memory. Instead of training after each observed experience, we wait for some iterations before sampling a batch of experiences from the replay memory and then train the neural networks with mini-batches. This batch training affords more training stability, allows the reuse of past experiences, and speeds up the training step.

Finally, to make better use of available information, particularly requests that are already revealed at a given time step, we propose an extension to using only Q-values. A look-ahead rollout strategy is proposed. We present this method in Section \ref{sec:rollout} and the results of our Q-learning application with and without rollouts is presented in Section \ref{sec:simulation_study}.

\subsection{Look-ahead rollout strategy}
\label{sec:rollout}

Monte Carlo Tree Search (MCTS) methods are exploration techniques to improve policies, based on previous work on Monte Carlo methods. \cite{coulom2007} first describes a tree search method applied to games. Several variants have been adapted for numerous applications, such as the famous AlphaGo Zero from \cite{Silver2017}. The general idea behind MCTS is to explore the action space by simulations following a tree structure, further exploring the promising areas of the tree.

For the RMFS real-time storage assignment, the objective of a look-ahead search is twofold. First, as in other applications, it helps the agent take better actions by exploring numerous possible scenarios. Second, and more specific to our application, the look-ahead can leverage information that is already revealed but not necessarily given in the state representation. This is the case of future orders. So far, in Section \ref{system_representation}, only information about the next order the robot needs to retrieve is included in the state features. However, at a given time step, while not all future orders have been revealed, some have. The look-ahead can then use \emph{only} those already-revealed orders (even if the real future sequence will differ) to take a better informed first action. 

\begin{figure}
\centering
\includegraphics[width=60mm]{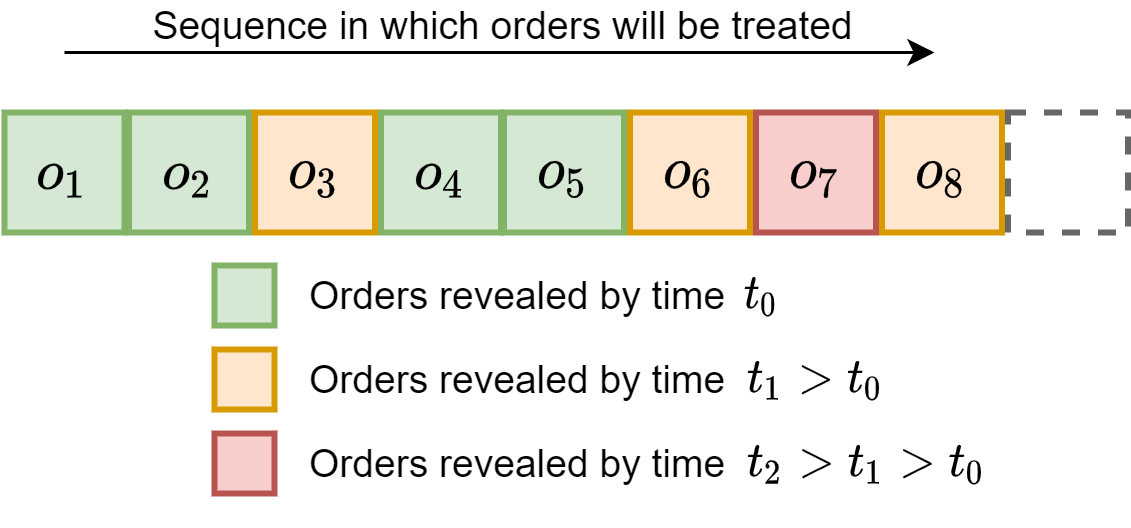}
\caption{Revealed orders at decision time}
\label{fig:mcts_revealed}
\end{figure}

Figure \ref{fig:mcts_revealed} illustrates this point. The orders $o_i$ are sorted according to the sequence in which they will actually be processed, but at time $t_0$ only the \emph{green} orders $\{o_1, \: o_2, \: o_4, \: o_5\}$ are known. The \emph{orange} orders $\{o_3, \: o_6, \: o_8\}$ will be revealed later at time $t_1>t_0$, for instance by the time order $o_2$ will have been processed. In this case, order $o_3$ is be inserted before $o_4$, for priority reasons. In this case, even if new orders will later be inserted into the list, at time $t_0$ only requests $\{o_1, \: o_2, \: o_4, \: o_5\}$ are used in the rollout, the only exceptions being potential other \emph{green} requests not depicted on Figure \ref{fig:mcts_revealed}.

\begin{figure}
\centering
\includegraphics[width=110mm]{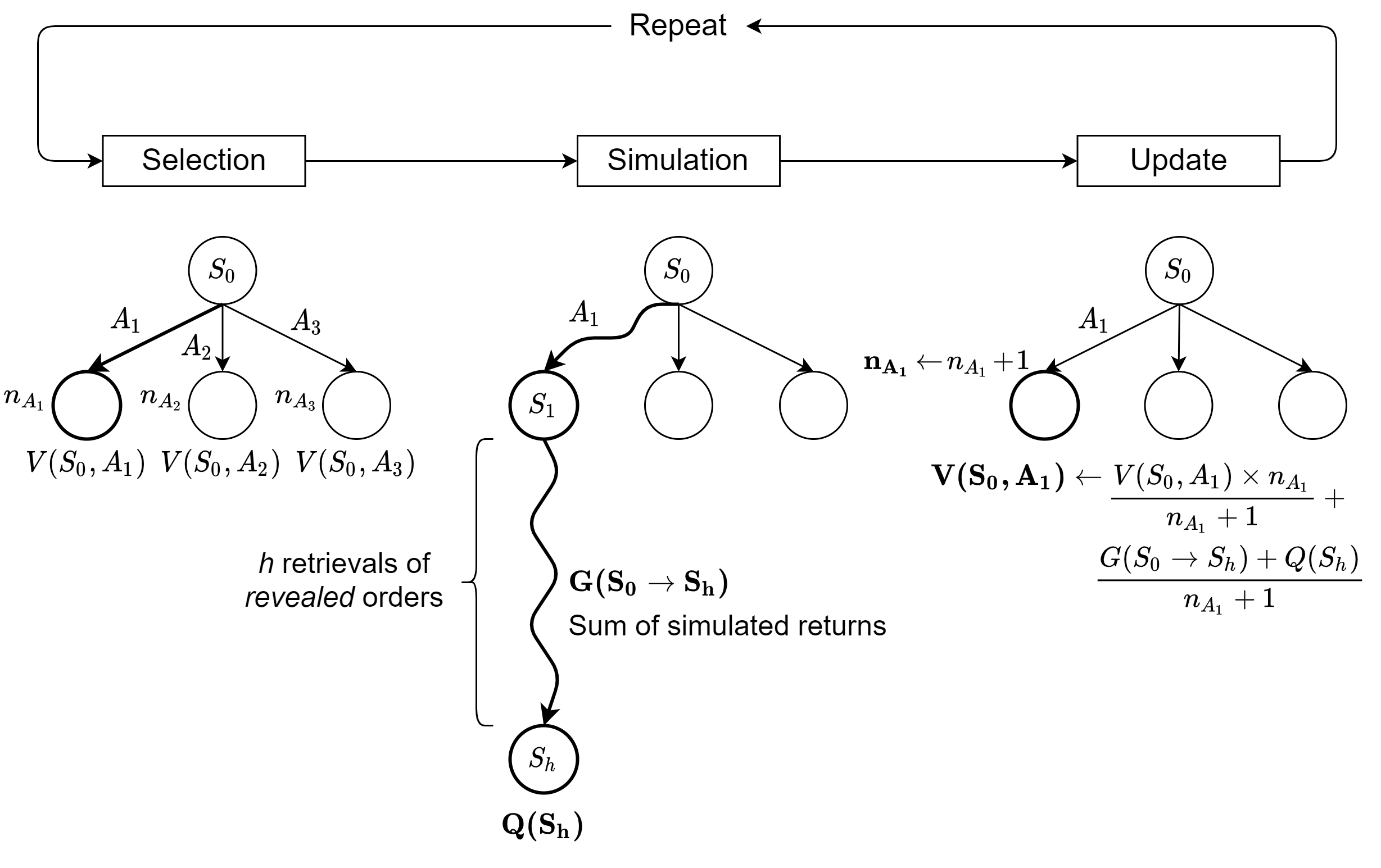}
\caption{Look-ahead rollout strategy}
\label{fig:mcts}
\end{figure}

\cite{Bertsekas2019ReinforcementControl} presents different approaches for MCTS implementation. One of these approaches is a single-step look-ahead with truncated rollout and cost approximation. Figure \ref{fig:mcts} presents the implementation of this type of rollout adapted to our problem. At first, all feasible actions at a given state $S_0$ (time $t_0$) are considered and their values $V(S_0, A_i)$ and visit counts $n_{A_i}$ are set to $0$. At each iteration, an action is randomly selected (action $A_1$ in Figure \ref{fig:mcts}) and a simulation of $h$ (horizon) consecutive storage tasks is run. This finite simulation uses orders revealed by time $t_0$. The action selection in the simulation uses the agent's current policy. The sum of reward along the simulation is computed, and the value of the last state $S_h$ is, by bootstrapping, set as the Q-value of the agent. Then, the value of the selected action is updated to its experienced average, and the visit counter is incremented. The process is repeated for a given number of iterations, and finally, the action presenting the best (min) value is selected. Note that in our specific usage of this look-ahead rollout strategy, each feasible action needs to be selected only once. Indeed, both the policy and the transitions within the rollout are deterministic because of the Q-learning agent's policy and because we do not simulate unknown future orders.

\section{Simulation study}
\label{sec:simulation_study}

In this section, we present a simulation study comparing our proposed methods to typical storage policies found in the literature. We start by defining the framework of the study, before presenting our results.

\subsection{Parameters}
\label{sec:parameters}

Kiva and Amazon have claimed travel speeds up to 3 to 4 miles per hour; however, to account for turns and congestion, we assume that the robots travel at a constant speed of 1.4 miles per hour (about 0.6 m/s). In our experiments, we assume an average picking time by the human operators of 8 sec and a loading and unloading time by the robots of 3 sec. Our study considers a relatively small storage area, with 36 storage locations (6 zones of size 2 by 3 shelves) for 34 shelves, 5 automated robots and one picking station.

Incoming orders in warehousing systems are typically skewed, with a small percentage of the items representing the vast majority of the orders and the others qualified as slow movers. \cite{Hausman1976OptimalSystems} propose modelling such a demand distribution with the definition of an \emph{ABC} curve: $G(x)=x^s \quad \text{for } 0 < s \leq 1$, which represents the proportion of cumulative demand (\%) of the first $x\%$ of the shelves. The smaller the skewness parameter $s$, the more skewed the demand distribution. We follow this suggestion by associating with every item type $i$ (equivalent to the corresponding shelf), at every discretised period, a Poisson distribution of average $\lambda _i = \left (\left ( \frac{i}{m} \right )^s - \left ( \frac{i-1}{m} \right )^s \right) \times \frac{n}{N}$, where $m$ is the number of item types, $n$ the total number of orders in the time horizon, and $N$ the number of discretised periods. Then, for every generated order $k$, a completion time $\delta_k$ is randomly drawn from a uniform distribution of interval $[1; \alpha \times T]$ where $\alpha$ defines the tightness of the completion times. Note that in the case where a drawn demand for a given item type is greater than one, the resulting orders are automatically grouped, corresponding to an opportunistic task of type 2 in \cite{Rimele2020RoboticApplications}. The values used here are: $m=34$, $n=30000$, $N=1440$ (corresponding to discretised periods every 1 min) for a time horizon of 24h, and $\alpha=0.4$. The values of the skewness parameters will be varied to study their impact on the different storage policies.

\begin{figure}
\centering
\includegraphics[width=80mm]{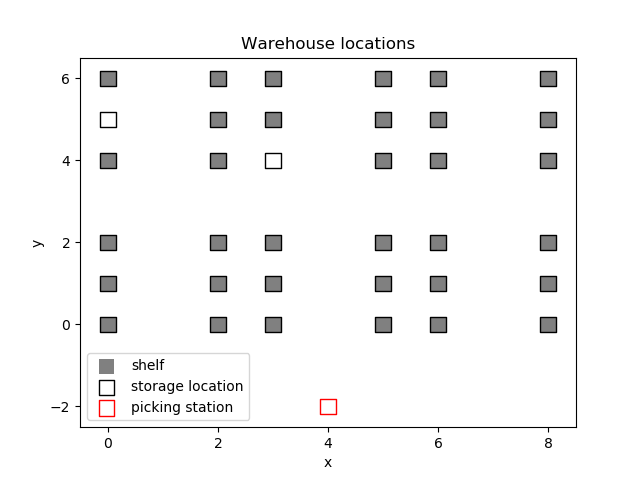}
\caption{Simulation study - Plan view of the storage area}
\label{fig:plan_view}
\end{figure}

After some trial and error, we set the parameter values for the RL agent as follows. Each training phase is conducted over 10000 episodes, corresponding to 24h of simulated incoming orders. Each of these training episodes uses its own generated instance of demand, while the comparisons between the baselines and the test of the RL policy are made on another shared instance. The exploration rate $\epsilon$ is kept constant at $0.1$. The neural networks are created using the Keras library. We use the mean square error loss function and the Adam optimiser for gradient descent. Both neural networks are feed-forward, fully connected with 3 inner layers of 32 neurons each and Relu activation functions. The output layer presents one neuron per zone of storage (6), and these neurons are linearly activated. The learning rate is set to $0.00025$, the replay memory has a capacity of $1000$, and batch training is performed every $100$ iterations. This batch consists of $256$ experiences sampled from the replay memory, and they are fitted to the model in mini-batches of $64$. The target network's weights are updated every $500$ iterations.

\subsection{Results}
\label{results}

This section presents results obtained by the proposed approach and the standard baselines. First, experiments were run using the method presented in Section \ref{sec:rl_method} on different instances corresponding to different values of demand skewness parameter $s$. Then, the rollout strategy presented in Section \ref{sec:rollout} was applied on both the Q-learning agent and the best-performing baseline, and gains that can be obtained by using additional available information were observed. The three baselines tested were the following: \emph{Random} storage policy, \emph{Class-based} and \emph{Shortest Leg}.

The \emph{Class-based} storage policy uses the concentric classes definition, depicted in Figure \ref{fig:class_zone_layouts} (left) and conventionally described in the literature. Since our method allocates shelves to zones instead of exact locations, the \emph{Shortest Leg} storage policy is implemented similarly. When a storage task needs to be performed and the location of the next retrieval is known, the exact potential location within each zone can be identified because of the arbitrary choice of selecting the location closest to the picking station. Knowing the candidate location from each zone (if any), the \emph{leg} distance can be computed as being the distance from the picking station to the storage location plus the distance from the storage location to the next retrieval. The selected zone corresponds to the smallest leg value. This policy results in greedily minimising the immediate cycle times.

\begin{figure}
\centering
\includegraphics[width=90mm]{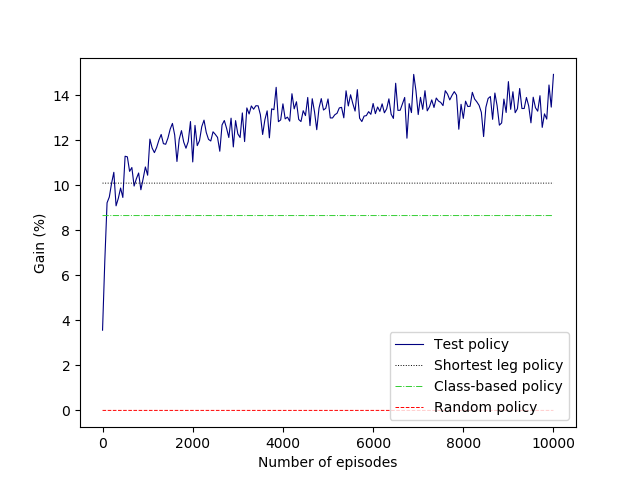}
\caption{Typical instance of a training curve (here for s=0.6)}
\label{fig:learning_curve}
\end{figure}

Table \ref{table:1} presents the average travelling times, $t(s)$, obtained from the three baselines and the RL agent, without rollouts, for several skewness parameter values ranging from $s=0.4$ (very skewed) to $s=1$ (flat distribution). For each storage policy, the corresponding performance gain, g(\%), compared to the \emph{Random} storage policy is also computed; this represents the percentage decrease in travelling time. As expected, the performance of the \emph{Class-based} policy increases when the skewness parameter value $s$ decreases. In concordance with the literature, we observe the good performance of the \emph{Shortest Leg} policy, which remains rather constant, between 10 and 11\% of improvement, depending on the skewness parameter values. In fact, it is only surpassed once by the \emph{Class-based} policy for the most skewed distribution $s=0.4$, but it offers the net advantage of being independent of the distribution pattern. The Q-learning agent consistently performs better than the baselines, between 3.21\% and 4.83\% of additional gain compared to the baselines relative to \emph{Random} storage. Similar to the SL policy, its performance remains rather independent of the skewness of the distribution. Figure \ref{fig:learning_curve} presents a typical training curve, here for parameter $s=0.6$. The \emph{x} axis shows the number of training episodes ranging from 0 to 10000 and \emph{y} corresponds to the relative gain in travelling time compared to the \emph{Random} storage policy. The horizontal dashed lines in the graph represent the performance of the three baseline policies. The test policy curve represents the performance of the RL policy in training, which is tested every 100 training episodes. With our training settings, we notice the non-linear evolution of the test policy curve until it reaches some plateau value around 5000 episodes.

\begin{table}
\centering
\small
\begin{tabular}{ |m{1.0cm}||m{0.9cm}|m{0.5cm}|m{0.9cm}|m{0.9cm}|m{0.9cm}|m{0.9cm}|m{0.9cm}|m{0.9cm}|  }
 \hline
 \multirow{2}{*}{\shortstack{Storage \\ policy}} & \multicolumn{2}{c|}{Random} & \multicolumn{2}{c|}{Class-based} & \multicolumn{2}{c|}{SL} & \multicolumn{2}{c|}{Agent} \\
 
            & t(s)       & g(\%)                       & t(s)       & g(\%)     & t(s)       & g(\%)   & t(s)   & g(\%) \\
\hline
 s=0.4      & 37.97      & \multicolumn{1}{c|}{-}      & 33.41      & 12.01      & 33.61      & 11.48   & 32.19  & 15.22 \\
 s=0.5      & 38.62      & \multicolumn{1}{c|}{-}      & 34.65      & 10.27      & 34.18      & 11.48   & 32.90  & 14.79 \\
 s=0.6      & 38.69      & \multicolumn{1}{c|}{-}      & 35.34      & 8.65      & 34.79      & 10.09   & 32.92  & 14.92 \\
 s=0.7      & 38.73      & \multicolumn{1}{c|}{-}      & 36.01      & 7.02      & 34.70      & 10.39   & 32.91  & 15.02 \\
 s=0.8      & 38.94      & \multicolumn{1}{c|}{-}      & 36.48      & 6.30      & 34.76      & 10.73   & 33.04  & 15.15 \\
 s=0.9      & 38.98      & \multicolumn{1}{c|}{-}      & 37.41      & 4.03      & 34.53      & 11.42   & 33.03  & 15.27 \\
 s=1.0      & 39.06      & \multicolumn{1}{c|}{-}      & 38.26      & 2.05      & 34.76      & 11.01   & 33.4  & 14.48 \\
 \hline
\end{tabular}
\caption{Average travelling times t(s) and performance gains g(\%) depending on distribution skewness parameter value s - without rollouts}
\label{table:1}
\end{table}

Over a second phase, we retained the overall best performing-baseline (SL) and our Q-learning agent to enhance them using the look-ahead rollout strategy presented in Section \ref{sec:rollout}. Again, we tested the resulting policies on different skewness parameter $s$ values, but we also considered different horizon values $h$, which define how far ahead the rollout is run. The only difference in the application of the rollout strategy to the SL baseline compared to the Q-learning agent is the absence of a Q-value at the last step of the rollout. In this case, the SL policy is run over the next $h$ already-revealed orders, and the average of the resulting cycle times is computed, for each first feasible action. Table \ref{table:2} presents the obtained results in terms of performance gains (\%) compared to \emph{Random} storage. First, we notice that for all parameter values, the performance is significantly improved. We can see that for both policies, the performance varies depending on the horizon of the rollout that is considered. As a general statement, it appears that horizon values 20 and 30 give the best results, with some exceptions. The fact that the performance increases with the horizon value before decreasing again is not too surprising. While using revealed orders in the rollouts is beneficial to making better-informed decisions, the longer the rollout, the more likely it is that new, unknown orders will be later inserted between the orders considered in the rollout. When this phenomenon starts to occur too often, it deteriorates the insight of the rollout and results in worse policies. The values in bold in the table designate the best results for each parameter $s$ for both policies. For the SL storage policy, the application of the rollout strategy increases its performance gains between 5.26\% ($s=1.0$) and 7.11\% ($s=0.7$). For the Q-learning agent, the performance gains range between 2.78\% ($s=0.8$) and 4.54\% ($s=0.4$). It appears that using rollouts benefits strongly skewed distributions slightly more than less skewed ones, which is most likely explained by more opportunities associated with fast and slow-moving shelves. Also, the rollout strategy benefits the SL policy more so than the agent. The differences in performance gains between the two policies with rollouts are now 2.29, 0.83, 1.47, 1.07, 1.55, 2.2 and 1.29\% for increasing $s$ values. Using rollouts gives a policy insights about future behaviours by simulating possible trajectories. While the Q-learning agent still benefits from this insight, with the Bellman equation it is, by design, already taking into account future impacts of a decision. On the other hand, the SL policy acts perfectly greedily; it has more to gain by looking into possible future outcomes. While the combination of the Q-learning agent with rollouts still performs best with an average of 1.53\% of performance gain compared to SL with rollouts, this assessment has interesting practical value because of the relative simplicity of implementation of SL with rollouts and its highly competitive performance. Finally, when comparing the results of the complete method with the best baseline commonly used from the literature (SL without rollouts), an average performance gain of 7.58\% is found.

\begin{table}
\centering
\footnotesize
\begin{tabular}{ |m{1cm}||m{0.6cm}|m{0.6cm}|m{0.65cm}|m{0.65cm}|m{0.65cm}|m{0.6cm}|m{0.6cm}|m{0.65cm}|m{0.65cm}|m{0.65cm}|  }
 \hline
  & \multicolumn{5}{c|}{SL + rollouts} & \multicolumn{5}{c|}{Agent + rollouts} \\
 
 h value    & \multicolumn{1}{c|}{5}     & \multicolumn{1}{c|}{10}    & \multicolumn{1}{c|}{20}    & \multicolumn{1}{c|}{30}    & \multicolumn{1}{c|}{40}    & \multicolumn{1}{c|}{5}     & \multicolumn{1}{c|}{10}    & \multicolumn{1}{c|}{20}    & \multicolumn{1}{c|}{30}    & \multicolumn{1}{c|}{40} \\
\hline
 s=0.4      & 14.01 & 16.83 & \textbf{17.47} & 17.07 & 17.22 & 16.26 & 17.01 & \textbf{19.76} & 18.84 & 18.80 \\
 s=0.5      & 14.91 & 16.03 & 17.57 & \textbf{17.65} & 16.00 & 15.17 & 17.14 & \textbf{18.48} & 18.41 & 17.89 \\
 s=0.6      & 14.97 & 15.47 & 16.06 & \textbf{16.92} & 16.05 & 15.71 & 17.45 & \textbf{18.39} & 17.99 & 18.16 \\
 s=0.7      & 14.41 & 15.45 & \textbf{17.50} & 16.22 & 15.85 & 14.96 & 16.82 & 18.13 & \textbf{18.57} & 18.11 \\
 s=0.8      & 14.23 & 15.51 & \textbf{16.38} & 15.42 & 15.46 & 15.13 & 16.32 & 17.23 & \textbf{17.93} & 16.96 \\
 s=0.9      & 14.22 & 15.49 & 16.34 & \textbf{16.79} & 16.52 & 15.11 & 17.04 & 17.83 & 18.95 & \textbf{18.99} \\
 s=1.0      & 15.15 & 15.60 & \textbf{16.27} & 15.64 & 15.61 & 15.27 & 16.81 & 17.34 & \textbf{17.56} & 17.17 \\
 \hline
\end{tabular}
\caption{Performance gains g(\%) depending on distribution skewness parameter value s and horizon h - with look-ahead rollouts}
\label{table:2}
\end{table}

\section{Conclusions}
\label{conclusions}

This work tackles the problem of dynamically allocating shelves to storage locations within a Robotic Mobile Fulfillment System. Because of the very nature of the e-commerce market, new orders need to be considered, if not fulfilled, as soon as they are revealed, which limits opportunities for batching. Typical methods from the literature correspond to decision rules which rely either on the distinct turnover rates of the containers or on minimising immediate cycles greedily. After making several assumptions about other operational decision rules in such warehouses, as well as the physical characteristics of the warehouse, we propose defining a storage policy using a POMDP and a Q-learning agent from Reinforcement Learning, to minimise the average travelling cycle time. This agent learns from repeated experiences which storage decision should be made based on a set of features representing the current state of the warehouse. Using zone-based storage allocation, we compared our method to decision rule baselines, including the \emph{Class-based} storage policy and \emph{Shortest Leg} storage policy. A performance gain between 3.5\% and 5\% higher than the best baseline relative to \emph{Random} storage was observed. In a second phase, the objective was to leverage additional information regarding orders that are already revealed at a given time step. While new orders will later appear and be inserted between those revealed orders, the latter can still be used in a look-ahead rollout strategy applied at the decision-making step. This rollout simulates finite horizon trajectories using the storage policy at hand for action selection within the trajectory, as well as a Q-value estimation at the last step, acting as an estimation of future expected objective value. Such a rollout strategy has the benefit of being applicable to any storage policy, which allowed us to make fair comparisons between the Q-learning agent and the \emph{Shortest Leg} policy, both using rollouts. After selecting the proper horizon length, using look-ahead rollouts increased the performance gains of the SL and Q-learning policies by about 6 and 4\%, respectively. Even though the Q-learning with rollouts policy gave the best results, the higher positive impact of rollouts on the SL policy makes it an interesting contender, due to its relative simplicity of application. Comparing the policy obtained from Q-learning with rollouts to the best baseline, we observed an average performance gain of 7.5\%, which, we think, demonstrates the potential for more research in this direction.

Future research may build on this work to extend it further. First, storage allocation could select exact open locations instead of using the notion of zones as described here. Instead of using single-step look-ahead rollouts, one could implement a complete Monte Carlo Tree Search algorithm to better explore action selection. The expensive computational cost of such an approach may limit online deployment but, coupled with a policy gradient agent instead of a Q-learning one, the tree search could be used extensively in the training phase only for fast future deployment. Another aspect worth looking into would be the waiting time of robots at the picking stations. In this work, we only consider the travelling time, which is a lower bound to the full cycle time. Aiming at minimising the travelling time along with implementing a distributed arrival of robots at the picking stations could result in higher throughput, but how to do that exactly appears to be quite challenging.

\section*{Acknowledgements}

We are grateful to the Mitacs Accelerate Program for providing funding for this project. We also wish to gratefully acknowledge the support and valuable insights of our industrial partners, JDA Software and Element AI.

%\section*{References}

\bibliography{mybibfile}

\begin{thebibliography}{39}
\expandafter\ifx\csname natexlab\endcsname\relax\def\natexlab#1{#1}\fi
\providecommand{\url}[1]{\texttt{#1}}
\providecommand{\href}[2]{#2}
\providecommand{\path}[1]{#1}
\providecommand{\DOIprefix}{doi:}
\providecommand{\ArXivprefix}{arXiv:}
\providecommand{\URLprefix}{URL: }
\providecommand{\Pubmedprefix}{pmid:}
\providecommand{\doi}[1]{\href{http://dx.doi.org/#1}{\path{#1}}}
\providecommand{\Pubmed}[1]{\href{pmid:#1}{\path{#1}}}
\providecommand{\bibinfo}[2]{#2}
\ifx\xfnm\relax \def\xfnm[#1]{\unskip,\space#1}\fi
%Type = Article
\bibitem[{Azadeh et~al.(2017)Azadeh, de~Koster \&
  Roy}]{Azadeh2017RobotizedOpportunities}
\bibinfo{author}{Azadeh, K.}, \bibinfo{author}{de~Koster, R.}, \&
  \bibinfo{author}{Roy, D.} (\bibinfo{year}{2017}).
\newblock \bibinfo{title}{{Robotized Warehouse Systems: Developments and
  Research Opportunities}}.
\newblock {\it \bibinfo{journal}{SSRN Electronic Journal}\/},  (pp.
  \bibinfo{pages}{1--55}). \DOIprefix\doi{10.2139/ssrn.2977779}.
%Type = Article
\bibitem[{Banker(2018)}]{Banker2018NewAutomation}
\bibinfo{author}{Banker, S.} (\bibinfo{year}{2018}).
\newblock \bibinfo{title}{{New Solution Changes the Rules of Warehouse
  Automation}}.
\newblock {\it \bibinfo{journal}{Forbes Business}\/}, . \URLprefix
  \url{https://www.forbes.com/sites/stevebanker/2018/06/12/new-solution-changes-the-rules-of-warehouse-automation}.
%Type = Article
\bibitem[{van~den Berg \& Gademann(2000)}]{VandenBerg2000SimulationSystem}
\bibinfo{author}{van~den Berg, J.~P.}, \& \bibinfo{author}{Gademann, A.}
  (\bibinfo{year}{2000}).
\newblock \bibinfo{title}{{Simulation study of an automated storage/retrieval
  system}}.
\newblock {\it \bibinfo{journal}{International Journal of Production
  Research}\/},  {\it \bibinfo{volume}{38}\/}, \bibinfo{pages}{1339--1356}.
  \DOIprefix\doi{10.1080/002075400188889}.
%Type = Article
\bibitem[{Bertsekas(1976)}]{BERTSEKAS1976DynamicControl}
\bibinfo{author}{Bertsekas, D.} (\bibinfo{year}{1976}).
\newblock \bibinfo{title}{{Dynamic Programming and Stochastic Control}}.
\newblock {\it \bibinfo{journal}{Mathematics in Science and Engineering}\/},
  {\it \bibinfo{volume}{125}\/}, \bibinfo{pages}{222--293}.
%Type = Book
\bibitem[{Bertsekas(2019)}]{Bertsekas2019ReinforcementControl}
\bibinfo{author}{Bertsekas, D.} (\bibinfo{year}{2019}).
\newblock {\it \bibinfo{title}{{Reinforcement learning and optimal
  control}}\/}.
\newblock \bibinfo{publisher}{Athena scientific}.
%Type = Article
\bibitem[{Boysen et~al.(2017)Boysen, Briskorn \&
  Emde}]{Boysen2017Parts-to-pickerEnvironment}
\bibinfo{author}{Boysen, N.}, \bibinfo{author}{Briskorn, D.}, \&
  \bibinfo{author}{Emde, S.} (\bibinfo{year}{2017}).
\newblock \bibinfo{title}{{Parts-to-picker based order processing in a
  rack-moving mobile robots environment}}.
\newblock {\it \bibinfo{journal}{European Journal of Operational Research}\/},
  {\it \bibinfo{volume}{262}\/}, \bibinfo{pages}{550--562}.
  \DOIprefix\doi{10.1016/j.ejor.2017.03.053}.
%Type = Article
\bibitem[{Boysen et~al.(2018)Boysen, Fedtke \&
  Weidinger}]{Boysen2018OptimizingProblem}
\bibinfo{author}{Boysen, N.}, \bibinfo{author}{Fedtke, S.}, \&
  \bibinfo{author}{Weidinger, F.} (\bibinfo{year}{2018}).
\newblock \bibinfo{title}{{Optimizing automated sorting in warehouses: The
  minimum order spread sequencing problem}}.
\newblock {\it \bibinfo{journal}{European Journal of Operational Research}\/},
  {\it \bibinfo{volume}{270}\/}, \bibinfo{pages}{386--400}.
  \DOIprefix\doi{10.1016/j.ejor.2018.03.026}.
%Type = Article
\bibitem[{Boysen et~al.(2019{\natexlab{a}})Boysen, de~Koster \&
  Weidinger}]{Boysen2019WarehousingSurvey}
\bibinfo{author}{Boysen, N.}, \bibinfo{author}{de~Koster, R.}, \&
  \bibinfo{author}{Weidinger, F.} (\bibinfo{year}{2019}{\natexlab{a}}).
\newblock \bibinfo{title}{{Warehousing in the e-commerce era: A survey}}.
\newblock {\it \bibinfo{journal}{European Journal of Operational Research}\/},
  {\it \bibinfo{volume}{277}\/}, \bibinfo{pages}{396--411}.
  \DOIprefix\doi{10.1016/j.ejor.2018.08.023}.
%Type = Article
\bibitem[{Boysen et~al.(2019{\natexlab{b}})Boysen, Stephan \&
  Weidinger}]{Boysen2019ManualProblem}
\bibinfo{author}{Boysen, N.}, \bibinfo{author}{Stephan, K.}, \&
  \bibinfo{author}{Weidinger, F.} (\bibinfo{year}{2019}{\natexlab{b}}).
\newblock \bibinfo{title}{{Manual order consolidation with put walls: the
  batched order bin sequencing problem}}.
\newblock {\it \bibinfo{journal}{EURO Journal on Transportation and
  Logistics}\/},  {\it \bibinfo{volume}{8}\/}, \bibinfo{pages}{169--193}.
  \DOIprefix\doi{10.1007/s13676-018-0116-0}.
%Type = Article
\bibitem[{Bozer \& White(1984)}]{Bozer1984Travel-TimeSystems}
\bibinfo{author}{Bozer, Y.~A.}, \& \bibinfo{author}{White, J.~A.}
  (\bibinfo{year}{1984}).
\newblock \bibinfo{title}{{Travel-Time Models for Automated Storage/Retrieval
  Systems}}.
\newblock {\it \bibinfo{journal}{IIE Transactions}\/},  {\it
  \bibinfo{volume}{16}\/}, \bibinfo{pages}{329--338}.
  \DOIprefix\doi{10.1080/07408178408975252}.
%Type = Inproceedings
\bibitem[{Coulom(2007)}]{coulom2007}
\bibinfo{author}{Coulom, R.} (\bibinfo{year}{2007}).
\newblock \bibinfo{title}{{Efficient Selectivity and Backup Operators in
  Monte-Carlo Tree Search}}.
\newblock In \bibinfo{editor}{H.~J. van~den Herik},
  \bibinfo{editor}{P.~Ciancarini}, \& \bibinfo{editor}{H.~H. L. M.~J. Donkers}
  (Eds.), {\it \bibinfo{booktitle}{Computers and Games}\/} (pp.
  \bibinfo{pages}{72--83}).
\newblock \bibinfo{address}{Berlin, Heidelberg}: \bibinfo{publisher}{Springer
  Berlin Heidelberg}.
%Type = Article
\bibitem[{D'Andrea \& Wurman(2008)}]{DAndrea2008FutureFacilities}
\bibinfo{author}{D'Andrea, R.}, \& \bibinfo{author}{Wurman, P.}
  (\bibinfo{year}{2008}).
\newblock \bibinfo{title}{{Future challenges of coordinating hundreds of
  autonomous vehicles in distribution facilities}}.
\newblock {\it \bibinfo{journal}{2008 IEEE International Conference on
  Technologies for Practical Robot Applications, TePRA}\/},  (pp.
  \bibinfo{pages}{80--83}). \DOIprefix\doi{10.1109/TEPRA.2008.4686677}.
%Type = Article
\bibitem[{Davarzani \& Norrman(2015)}]{Davarzani2015}
\bibinfo{author}{Davarzani, H.}, \& \bibinfo{author}{Norrman, A.}
  (\bibinfo{year}{2015}).
\newblock \bibinfo{title}{{Toward a relevant agenda for warehousing research:
  literature review and practitioners’ input}}.
\newblock {\it \bibinfo{journal}{Logistics Research}\/},  {\it
  \bibinfo{volume}{8}\/}. \DOIprefix\doi{10.1007/s12159-014-0120-1}.
%Type = Article
\bibitem[{Enright \& Wurman(2011)}]{Enright2011OptimizationSystems}
\bibinfo{author}{Enright, J.~J.}, \& \bibinfo{author}{Wurman, P.~R.}
  (\bibinfo{year}{2011}).
\newblock \bibinfo{title}{{Optimization and coordinated autonomy in mobile
  fulfillment systems}}.
\newblock {\it \bibinfo{journal}{AAAI Workshop - Technical Report}\/},  {\it
  \bibinfo{volume}{WS-11-09}\/}, \bibinfo{pages}{33--38}.
%Type = Article
\bibitem[{Gagliardi et~al.(2012)Gagliardi, Renaud \&
  Ruiz}]{Gagliardi2012OnSystems}
\bibinfo{author}{Gagliardi, J.-P.}, \bibinfo{author}{Renaud, J.}, \&
  \bibinfo{author}{Ruiz, A.} (\bibinfo{year}{2012}).
\newblock \bibinfo{title}{{On storage assignment policies for unit-load
  automated storage and retrieval systems}}.
\newblock {\it \bibinfo{journal}{International Journal of Production
  Research}\/},  {\it \bibinfo{volume}{50}\/}, \bibinfo{pages}{879--892}.
  \DOIprefix\doi{10.1080/00207543.2010.543939}.
%Type = Article
\bibitem[{Graves et~al.(1977)Graves, Hausman \&
  Schwarz}]{Graves1977Storage-RetrievalSystems}
\bibinfo{author}{Graves, S.~C.}, \bibinfo{author}{Hausman, W.~H.}, \&
  \bibinfo{author}{Schwarz, L.~B.} (\bibinfo{year}{1977}).
\newblock \bibinfo{title}{{Storage-Retrieval Interleaving in Automatic
  Warehousing Systems}}.
\newblock {\it \bibinfo{journal}{Management Science}\/},  {\it
  \bibinfo{volume}{23}\/}, \bibinfo{pages}{935--945}.
  \DOIprefix\doi{10.1287/mnsc.23.9.935}.
%Type = Article
\bibitem[{Gu et~al.(2007)Gu, Goetschalckx \& McGinnis}]{Gu2007}
\bibinfo{author}{Gu, J.}, \bibinfo{author}{Goetschalckx, M.}, \&
  \bibinfo{author}{McGinnis, L.~F.} (\bibinfo{year}{2007}).
\newblock \bibinfo{title}{{Research on warehouse operation: A comprehensive
  review}}.
\newblock {\it \bibinfo{journal}{European Journal of Operational Research}\/},
  {\it \bibinfo{volume}{177}\/}, \bibinfo{pages}{1--21}.
  \DOIprefix\doi{10.1016/j.ejor.2006.02.025}.
%Type = Article
\bibitem[{van Hasselt(2010)}]{VanHasselt2010DoubleQ-learning}
\bibinfo{author}{van Hasselt, H.} (\bibinfo{year}{2010}).
\newblock \bibinfo{title}{{Double Q-learning}}.
\newblock {\it \bibinfo{journal}{Advances in Neural Information Processing
  Systems 23: 24th Annual Conference on Neural Information Processing Systems
  2010, NIPS 2010}\/},  (pp. \bibinfo{pages}{1--9}).
%Type = Article
\bibitem[{van Hasselt et~al.(2016)van Hasselt, Guez \&
  Silver}]{vanHasselt2016DeepQ-Learning}
\bibinfo{author}{van Hasselt, H.}, \bibinfo{author}{Guez, A.}, \&
  \bibinfo{author}{Silver, D.} (\bibinfo{year}{2016}).
\newblock \bibinfo{title}{{Deep Reinforcement Learning with Double
  Q-Learning}}.
\newblock {\it \bibinfo{journal}{Proceedings of the 30th AAAI Conference on
  Artificial Intelligence (AAAI-16)}\/},  (pp. \bibinfo{pages}{2094--2100}).
%Type = Article
\bibitem[{Hausman et~al.(1976)Hausman, Schwarz \&
  Graves}]{Hausman1976OptimalSystems}
\bibinfo{author}{Hausman, W.~H.}, \bibinfo{author}{Schwarz, L.~B.}, \&
  \bibinfo{author}{Graves, S.~C.} (\bibinfo{year}{1976}).
\newblock \bibinfo{title}{{Optimal Storage Assignment in Automatic Warehousing
  Systems}}.
\newblock {\it \bibinfo{journal}{Management Science}\/},  {\it
  \bibinfo{volume}{22}\/}, \bibinfo{pages}{629--638}.
  \DOIprefix\doi{10.1287/mnsc.22.6.629}.
%Type = Inproceedings
\bibitem[{Hessel et~al.(2018)Hessel, Modayil, Van~Hasselt, Schaul, Ostrovski,
  Dabney, Horgan, Piot, Azar \& Silver}]{Hessel2018Rainbow:Learning}
\bibinfo{author}{Hessel, M.}, \bibinfo{author}{Modayil, J.},
  \bibinfo{author}{Van~Hasselt, H.}, \bibinfo{author}{Schaul, T.},
  \bibinfo{author}{Ostrovski, G.}, \bibinfo{author}{Dabney, W.},
  \bibinfo{author}{Horgan, D.}, \bibinfo{author}{Piot, B.},
  \bibinfo{author}{Azar, M.}, \& \bibinfo{author}{Silver, D.}
  (\bibinfo{year}{2018}).
\newblock \bibinfo{title}{{Rainbow: Combining improvements in deep
  reinforcement learning}}.
\newblock In {\it \bibinfo{booktitle}{32nd AAAI Conference on Artificial
  Intelligence, AAAI 2018}\/} (pp. \bibinfo{pages}{3215--3222}).
%Type = Article
\bibitem[{Kirks et~al.(2012)Kirks, Stenzel, Kamagaew \&
  en~Hompel}]{Kirks2012CellularSystems}
\bibinfo{author}{Kirks, T.}, \bibinfo{author}{Stenzel, J.},
  \bibinfo{author}{Kamagaew, A.}, \& \bibinfo{author}{en~Hompel, M.}
  (\bibinfo{year}{2012}).
\newblock \bibinfo{title}{{Cellular Transport Vehicles for Flexible and
  Changeable Facility Logistics Systems}}.
\newblock {\it \bibinfo{journal}{Logistics Journal}\/},  {\it
  \bibinfo{volume}{2192(9084)}\/}.
%Type = Article
\bibitem[{Lamballais et~al.(2017)Lamballais, Roy \&
  de~Koster}]{Lamballais2017EstimatingSystem}
\bibinfo{author}{Lamballais, T.}, \bibinfo{author}{Roy, D.}, \&
  \bibinfo{author}{de~Koster, R.} (\bibinfo{year}{2017}).
\newblock \bibinfo{title}{{Estimating performance in a Robotic Mobile
  Fulfillment System}}.
\newblock {\it \bibinfo{journal}{European Journal of Operational Research}\/},
  {\it \bibinfo{volume}{256}\/}, \bibinfo{pages}{976--990}.
  \DOIprefix\doi{10.1016/j.ejor.2016.06.063}.
%Type = Article
\bibitem[{Merschformann et~al.(2019)Merschformann, Lamballais, de~Koster \&
  Suhl}]{Merschformann2019DecisionSystems}
\bibinfo{author}{Merschformann, M.}, \bibinfo{author}{Lamballais, T.},
  \bibinfo{author}{de~Koster, R.}, \& \bibinfo{author}{Suhl, L.}
  (\bibinfo{year}{2019}).
\newblock \bibinfo{title}{{Decision rules for robotic mobile fulfillment
  systems}}.
\newblock {\it \bibinfo{journal}{Operations Research Perspectives}\/},  {\it
  \bibinfo{volume}{6}\/}, \bibinfo{pages}{100128}.
  \DOIprefix\doi{10.1016/j.orp.2019.100128}.
%Type = Article
\bibitem[{Merschformann et~al.(2018)Merschformann, Xie \&
  Li}]{Merschformann2018RAWSim-O:Systems}
\bibinfo{author}{Merschformann, M.}, \bibinfo{author}{Xie, L.}, \&
  \bibinfo{author}{Li, H.} (\bibinfo{year}{2018}).
\newblock \bibinfo{title}{{RAWSim-O: A simulation framework for robotic mobile
  fulfillment systems}}.
\newblock {\it \bibinfo{journal}{Logistics Research}\/},  {\it
  \bibinfo{volume}{11}\/}, \bibinfo{pages}{1--11}.
  \DOIprefix\doi{10.23773/2018{\_}8}.
%Type = Article
\bibitem[{Mnih et~al.(2015)Mnih, Kavukcuoglu, Silver, Rusu, Veness, Bellemare,
  Graves, Riedmiller, Fidjeland, Ostrovski, Petersen, Beattie, Sadik,
  Antonoglou, King, Kumaran, Wierstra, Legg \& Hassabis}]{Mnih2015}
\bibinfo{author}{Mnih, V.}, \bibinfo{author}{Kavukcuoglu, K.},
  \bibinfo{author}{Silver, D.}, \bibinfo{author}{Rusu, A.~A.},
  \bibinfo{author}{Veness, J.}, \bibinfo{author}{Bellemare, M.~G.},
  \bibinfo{author}{Graves, A.}, \bibinfo{author}{Riedmiller, M.},
  \bibinfo{author}{Fidjeland, A.~K.}, \bibinfo{author}{Ostrovski, G.},
  \bibinfo{author}{Petersen, S.}, \bibinfo{author}{Beattie, C.},
  \bibinfo{author}{Sadik, A.}, \bibinfo{author}{Antonoglou, I.},
  \bibinfo{author}{King, H.}, \bibinfo{author}{Kumaran, D.},
  \bibinfo{author}{Wierstra, D.}, \bibinfo{author}{Legg, S.}, \&
  \bibinfo{author}{Hassabis, D.} (\bibinfo{year}{2015}).
\newblock \bibinfo{title}{{Human-level control through deep reinforcement
  learning}}.
\newblock {\it \bibinfo{journal}{Nature}\/},  {\it \bibinfo{volume}{518}\/},
  \bibinfo{pages}{529--533}. \DOIprefix\doi{10.1038/nature14236}.
%Type = Incollection
\bibitem[{Park(2012)}]{Park2012OrderModels}
\bibinfo{author}{Park, B.~C.} (\bibinfo{year}{2012}).
\newblock \bibinfo{title}{{Order Picking: Issues, Systems and Models}}.
\newblock In {\it \bibinfo{booktitle}{Warehousing in the Global Supply
  Chain}\/} (pp. \bibinfo{pages}{1--30}).
\newblock \bibinfo{address}{London}: \bibinfo{publisher}{Springer London}.
\newblock \DOIprefix\doi{10.1007/978-1-4471-2274-6{\_}1}.
%Type = Article
\bibitem[{Ramanathan(2010)}]{Ramanathan2010TheE-commerce}
\bibinfo{author}{Ramanathan, R.} (\bibinfo{year}{2010}).
\newblock \bibinfo{title}{{The moderating roles of risk and efficiency on the
  relationship between logistics performance and customer loyalty in
  e-commerce}}.
\newblock {\it \bibinfo{journal}{Transportation Research Part E: Logistics and
  Transportation Review}\/},  {\it \bibinfo{volume}{46}\/},
  \bibinfo{pages}{950--962}. \DOIprefix\doi{10.1016/j.tre.2010.02.002}.
%Type = Article
\bibitem[{Rim{\'{e}}l{\'{e}} et~al.(2020)Rim{\'{e}}l{\'{e}}, Gamache, Gendreau,
  Grangier \& Rousseau}]{Rimele2020RoboticApplications}
\bibinfo{author}{Rim{\'{e}}l{\'{e}}, A.}, \bibinfo{author}{Gamache, M.},
  \bibinfo{author}{Gendreau, M.}, \bibinfo{author}{Grangier, P.}, \&
  \bibinfo{author}{Rousseau, L.-M.} (\bibinfo{year}{2020}).
\newblock \bibinfo{title}{{Robotic Mobile Fulfillment Systems: a Mathematical
  Modelling Framework for E-commerce Applications}}.
\newblock {\it \bibinfo{journal}{Cahiers du Cirrelt}\/},  {\it
  \bibinfo{volume}{CIRRELT-2020-42}\/}.
%Type = Article
\bibitem[{Roodbergen \& Vis(2009)}]{Roodbergen2009}
\bibinfo{author}{Roodbergen, K.~J.}, \& \bibinfo{author}{Vis, I. F.~A.}
  (\bibinfo{year}{2009}).
\newblock \bibinfo{title}{{A survey of literature on automated storage and
  retrieval systems}}.
\newblock {\it \bibinfo{journal}{European Journal of Operational Research}\/},
  {\it \bibinfo{volume}{194}\/}, \bibinfo{pages}{343--362}.
  \DOIprefix\doi{10.1016/j.ejor.2008.01.038}.
%Type = Article
\bibitem[{Roy et~al.(2019)Roy, Nigam, de~Koster, Adan \& Resing}]{Roy2019}
\bibinfo{author}{Roy, D.}, \bibinfo{author}{Nigam, S.},
  \bibinfo{author}{de~Koster, R.}, \bibinfo{author}{Adan, I.}, \&
  \bibinfo{author}{Resing, J.} (\bibinfo{year}{2019}).
\newblock \bibinfo{title}{{Robot-storage zone assignment strategies in mobile
  fulfillment systems}}.
\newblock {\it \bibinfo{journal}{Transportation Research Part E: Logistics and
  Transportation Review}\/},  {\it \bibinfo{volume}{122}\/},
  \bibinfo{pages}{119--142}. \DOIprefix\doi{10.1016/j.tre.2018.11.005}.
%Type = Article
\bibitem[{Schwarz et~al.(1978)Schwarz, Graves \&
  Hausman}]{Schwarz1978SchedulingResults}
\bibinfo{author}{Schwarz, L.~B.}, \bibinfo{author}{Graves, S.~C.}, \&
  \bibinfo{author}{Hausman, W.~H.} (\bibinfo{year}{1978}).
\newblock \bibinfo{title}{{Scheduling Policies for Automatic Warehousing
  Systems: Simulation Results}}.
\newblock {\it \bibinfo{journal}{A I I E Transactions}\/},  {\it
  \bibinfo{volume}{10}\/}, \bibinfo{pages}{260--270}.
  \DOIprefix\doi{10.1080/05695557808975213}.
%Type = Article
\bibitem[{Silver et~al.(2017)Silver, Schrittwieser, Simonyan, Antonoglou,
  Huang, Guez, Hubert, Baker, Lai, Bolton, Chen, Lillicrap, Hui, Sifre, {Van
  Den Driessche}, Graepel \& Hassabis}]{Silver2017}
\bibinfo{author}{Silver, D.}, \bibinfo{author}{Schrittwieser, J.},
  \bibinfo{author}{Simonyan, K.}, \bibinfo{author}{Antonoglou, I.},
  \bibinfo{author}{Huang, A.}, \bibinfo{author}{Guez, A.},
  \bibinfo{author}{Hubert, T.}, \bibinfo{author}{Baker, L.},
  \bibinfo{author}{Lai, M.}, \bibinfo{author}{Bolton, A.},
  \bibinfo{author}{Chen, Y.}, \bibinfo{author}{Lillicrap, T.},
  \bibinfo{author}{Hui, F.}, \bibinfo{author}{Sifre, L.}, \bibinfo{author}{{Van
  Den Driessche}, G.}, \bibinfo{author}{Graepel, T.}, \&
  \bibinfo{author}{Hassabis, D.} (\bibinfo{year}{2017}).
\newblock \bibinfo{title}{{Mastering the game of Go without human knowledge}}.
\newblock {\it \bibinfo{journal}{Nature}\/},  {\it \bibinfo{volume}{550}\/},
  \bibinfo{pages}{354--359}. \DOIprefix\doi{10.1038/nature24270}.
%Type = Book
\bibitem[{Sutton \& Barto(2018)}]{Sutton2018ReinforcementIntroduction}
\bibinfo{author}{Sutton, R.~S.}, \& \bibinfo{author}{Barto, A.~G.}
  (\bibinfo{year}{2018}).
\newblock {\it \bibinfo{title}{{Reinforcement Learning: An Introduction}}\/}.
\newblock \bibinfo{publisher}{MIT press}.
%Type = Phdthesis
\bibitem[{Watkins(1989)}]{Watkins1989LearningRewards}
\bibinfo{author}{Watkins, C.} (\bibinfo{year}{1989}).
\newblock {\it \bibinfo{title}{{Learning from delayed rewards}}\/}.
\newblock Ph.D. thesis King's College.
%Type = Article
\bibitem[{Weidinger et~al.(2018)Weidinger, Boysen \&
  Briskorn}]{Weidinger2018StorageWarehouses}
\bibinfo{author}{Weidinger, F.}, \bibinfo{author}{Boysen, N.}, \&
  \bibinfo{author}{Briskorn, D.} (\bibinfo{year}{2018}).
\newblock \bibinfo{title}{{Storage assignment with rack-moving mobile robots in
  KIVA warehouses}}.
\newblock {\it \bibinfo{journal}{Transportation Science}\/},  {\it
  \bibinfo{volume}{52}\/}, \bibinfo{pages}{1479--1495}.
  \DOIprefix\doi{10.1287/trsc.2018.0826}.
%Type = Article
\bibitem[{Wurman et~al.(2008)Wurman, D'Andrea \&
  Mountz}]{Wurman2008CoordinatingWarehouses}
\bibinfo{author}{Wurman, P.~R.}, \bibinfo{author}{D'Andrea, R.}, \&
  \bibinfo{author}{Mountz, M.} (\bibinfo{year}{2008}).
\newblock \bibinfo{title}{{Coordinating Hundreds of Cooperative, Autonomous
  Vehicles in Warehouses}}.
\newblock {\it \bibinfo{journal}{AI Magazine}\/},  {\it
  \bibinfo{volume}{29}\/}, \bibinfo{pages}{9--9}.
  \DOIprefix\doi{10.1609/AIMAG.V29I1.2082}.
%Type = Article
\bibitem[{Yaman et~al.(2012)Yaman, Karasan \& Kara}]{Yaman2012ReleaseDelivery}
\bibinfo{author}{Yaman, H.}, \bibinfo{author}{Karasan, O.~E.}, \&
  \bibinfo{author}{Kara, B.~Y.} (\bibinfo{year}{2012}).
\newblock \bibinfo{title}{{Release time scheduling and hub location for
  next-day delivery}}.
\newblock {\it \bibinfo{journal}{Operations Research}\/},  {\it
  \bibinfo{volume}{60}\/}, \bibinfo{pages}{906--917}.
  \DOIprefix\doi{10.1287/opre.1120.1065}.
%Type = Article
\bibitem[{Yuan et~al.(2018)Yuan, Graves \& Cezik}]{Yuan2018}
\bibinfo{author}{Yuan, R.}, \bibinfo{author}{Graves, S.~C.}, \&
  \bibinfo{author}{Cezik, T.} (\bibinfo{year}{2018}).
\newblock \bibinfo{title}{{Velocity-Based Storage Assignment in Semi-Automated
  Storage Systems}}.
\newblock {\it \bibinfo{journal}{Production and Operations Management}\/},
  {\it \bibinfo{volume}{28}\/}, \bibinfo{pages}{354--373}.
  \DOIprefix\doi{10.1111/poms.12925}.

\end{thebibliography}

\end{document}